\documentclass[12pt, reqno, twoside]{amsart}
\usepackage{amsmath, amsthm, amscd, amsfonts, amssymb, graphicx, color}
\usepackage{mathrsfs}
\usepackage[bookmarksnumbered, colorlinks, plainpages]{hyperref}
\hypersetup{colorlinks=true,linkcolor=red, anchorcolor=green, citecolor=cyan, urlcolor=red, filecolor=magenta, pdftoolbar=true}
\usepackage{enumitem}
\newcommand{\disp}{\displaystyle} 
\usepackage{enumitem}

\newtheorem{theorem}{Theorem}[section]
\newtheorem{lemma}[theorem]{Lemma}
\newtheorem{proposition}[theorem]{Proposition}
\newtheorem{corollary}[theorem]{Corollary}
\theoremstyle{definition}
\newtheorem{definition}[theorem]{Definition}

\theoremstyle{remark}
\newtheorem{remark}[theorem]{Remark}
\numberwithin{equation}{section}

\begin{document}

\setcounter{page}{1}

\title[Anisotropic conformable fractional Hardy type inequalities with weights]{Anisotropic Hardy type inequalities with weights and conformable fractional differential operators}

\author[A. Abolarinwa]{Abimbola Abolarinwa$^{*1}$}

\author[Y. O. Anthonio]{Yisa O. Anthonio$^{2}$}



\address{$^{*1}$Department of Mathematics,
University of Lagos, Akoka, Lagos State,
Nigeria}
\email{\textcolor[rgb]{0.00,0.00,0.84}{(Correspondence) a.abolarinwa1@gmail.com}}

\address{$^{2}$Department of Mathematics Science,  Lagos State University of Science and technology,  Ikorodu, Lagos State, Nigeria}

\date{November 2,  2024}

\begin{abstract}
By a systematic development of fundamental concepts of conformable calculus we establish conformable divergence theorem and Green's identities which we combine with some new anisotropic Picone type identities to derive  a generalized anisotropic Hardy type inequality with weights and conformable fractional differential operators. As a consequence,  several Hardy type inequalities and Heisenberg Pauli-Weyl uncertainty principles are obtained.
\end{abstract}
\maketitle
              

\section{Introduction}

Hardy type inequalities are one of important classical inequalities in mathematical analysis. G. H. Hardy introduced both discrete and continuous (integral) versions of these inequalities in 1920's.  Specifically, the discrete version was introduced in 1920 \cite{Hardy1} as follows:
\begin{align}\label{e11}
\sum_{m=1}^\infty \left(\frac{1}{m}\sum_{j=1}^m A_j\right)^p \le \left(\frac{p}{p-1}\right)^p \sum_{m=1}^\infty A^p_m
\end{align}
for $p>1$ and a nonnegative sequence of real numbers $\{A_m\}_{m=1}^\infty$. Here, both sides of \eqref{e11} must be finite, and equality holds if and only if $A_m=0$.  In 1925, Hardy \cite{Hardy2} proved also the integral version of \eqref{e11} for a nonnegative integrable function $f$ over $(0,x), x>0$ as 
\begin{align}\label{e12}
\int_0^\infty \left(\frac{1}{x}\int_0^x f(t)dt\right)^p dx \le \left(\frac{p}{p-1}\right)^p\int_0^\infty |f(x)|^p dx,
\end{align}
where $p>1$. Here $|f|^p$ is integrable and convergent over $(0,\infty)$ while the integral on the right hand side is finite, and equality holds if and only if $f(x)=0$.  In both cases \eqref{e11} and \eqref{e12} $(p/(p-1))^p$ is the best possible constant achievable.  For reference purpose, the history and developments of these inequalities are documented in the books \cite{ARSBook, KMPBook, OpicKBook, PSBook}. One can easily show that inequality \eqref{e12} is equivalent to 
\begin{align}\label{e13}
\int_0^\infty \frac{|f(x)|^p}{|x|^p} dx \le \left(\frac{p}{p-1}\right)^p\int_0^\infty \left|\frac{df(x)}{dx}\right|^p dx.
\end{align}

The integral version of Hardy inequalities has been extensively studied and used as a model example for the investigation of more general integral inequalities. An important generalization of Hardy inequalities is the multidimensional version which has numerous useful consequences in the theory of partial differential equations and differential geometry. For $f\in C^\infty_0(\mathbb{R}^n)$, $n\ge 2$, the $n$-dimensional version of  \eqref{e13} is given by 
\begin{align}\label{e14}
\int_\Omega \frac{|f(x)|^p}{|x|^p} dx    \le \left(\frac{p}{n-p}\right)^p \int_\Omega |\nabla f(x)|^p dx,
\end{align}
where $\Omega \subset \mathbb{R}^n$, $1\le p <n$ and $\nabla := (\partial_{x_1},  \partial_{x_2}, \cdots , \partial_{x_n})$ denotes the Euclidean gradient operator. Inequality \eqref{e14} can also be extended to the whole of $\mathbb{R}^n\setminus\{0\}$ for $p>n$, while the constant $(p/(n-p))^p$ is sharp but never achieved by a non-trivial function. For further discussion on \eqref{e14} and its variants, see \cite{AbA,ARY,Ab23,AbR,BV,DavH,KRS,RSBook,RSS,RV}. 
Hardy inequality of the form \eqref{e13} (or \eqref{e14}) is referred to as inequality of integer order, and has been studied extensively, not only in the Euclidean context,  but also in several other contexts like Heisenberg group, homogeneous groups, stratified Lie groups, Riemannian manifolds and so on \cite{Ab23, AbR, GL,KRS, RSBook,RSS}.

Recently, inequalities of fractional order have gained attention of researchers due to their applications in many scenarios involving non-locality.  Fractional derivatives of Riemann-Liouville and Caputo have been seriously engaged in this sense \cite{BD,JS,US}. The multi-dimensional fractional order Hardy inequality take the form \cite{Dyd,KRS2}
\begin{align}\label{e15}
\int_{\mathbb{R}^n} \frac{|f(x)|^p}{|x|^{\alpha p}} dx    \le C(n,\alpha,p) \iint_{\mathbb{R}^n \times {\mathbb{R}^n}} \frac{|f(x)-f(y)|^p}{|x-y|^{n+\alpha p}} dx dy
\end{align}
for $p>1$, $0<\alpha<1$ and $\alpha p<n$, $n>1$. The double integral on the right hand side is with respect to the fractional Gagliardo seminorm for a measurable function $f\in W^{\alpha,p}_0(\mathbb{R}^n)$.

There are shortcomings in the definitions of Riemann-Liouville and Caputo fractional derivatives in the sense that they lack certain vital properties associated with the classical derivative, such as properties that derivatives of constant should be zero, (Caputo fractional derivative of a constant is zero, though), composition rule,  Green's theorem and so on. These shortcomings limit their applicability to real life phenomena. To circumvent these challenges, Khalil et al \cite{Khal} introduced what is known as conformable fractional derivative of a function, thereby making it more flexible to accommodate many classical theorems of calculus, which in turn allows for extension of classical results to the fractional order set up. Further properties of conformable derivative are examined by authors in \cite{Abd,Ata}, while basic concepts of conformable fractional calculus are highlighted in Section \ref{sec2} of this paper.

The anisotropic conformable fractional differential operator is defined for  continuous $\alpha$-differentiable function $u$:\begin{align}\label{e16}
\sum_{k=1}^n\frac{\partial^\alpha }{\partial x_k^\alpha}\left(\left|\frac{\partial^\alpha u }{\partial x_k^\alpha}\right|^{p_k-2}\frac{\partial^\alpha u }{\partial x_k^\alpha} \right)
\end{align}
with $\alpha\in (0,1]$,  $p_k>1$, $k=1,\cdots,n$.  Setting $p_k=2$ and $p_k=p$ for all $k$ in \eqref{e16}, this operator reduces to  conformable Laplacian and the conformable pseudo-$p$-Laplacian,  respectively.  The anisotropic Laplacian plays crucial roles in several areas of mathematical theories and their applications in engineering and sciences.   For instance,  it reflects anisotropic characteristics of some reinforced materials \cite{Ta},  as well as explains fluid dynamics in the anisotropic media having different conductivities  in each direction \cite{Be}.  Models involving anisotropic Laplacian  arise also in image processing and computer vision \cite{Sa,Weic}.
Anisotropic Picone identities for classical gradient operator is proved in \cite{FC}  for differentiable functions $u\ge0$, $v>0$ in a domain of $\mathbb{R}^n$ and exponents $p_k>1$ as follows:
\begin{align}\label{e17}
\sum_{k=1}^n \left|\frac{\partial u}{\partial x_k}\right|^{p_k} - \sum_{k=1}^n p_k \left(\frac{u}{v}\right)^{p_k-1}  \left|\frac{\partial v}{\partial x_k}\right|^{p_k-2} \frac{\partial v}{\partial x_k}\frac{\partial u}{\partial x_k}
+ \sum_{k=1}^n  \left(\frac{u}{v}\right)^{p_k} \left|\frac{\partial v}{\partial x_k}\right|^{p_k} \ge 0
\end{align}
with equality if and only if $u=cv$ for some constant $c>0$. Picone type identities have proved to be effective tools in the study of existence and nonexistence of positive solutions to differential equations,  Sturmian comparison principle,  domain monotonicity,  Hardy's inequality, Caccioppoli inequality, e. t. c. (see \cite{AbP,AbB,AH,FC,Ja,Khel} and the references cited therein).

 This paper's main aim is therefore to establish anisotropic weighted fractional Hardy type inequalities by leveraging on the properties of conformable fractional derivative (see Theorem \ref{thm32}).  The generalized conformable fractional Hardy inequalities will come as a consequence of the anisotropic Picone identities.  So, we shall first prove some Picone type identities for anisotropic fractional gradient operator on a compatible domain (see Proposition \ref{prop31}).   One can see that the results of this paper are different from those obtained recently for single valued functions by authors in \cite{GSAZ1, GSAZ2}.  The method adopted here appears more general as we consider multidimensional functions and are able to take care of the boundary terms.  Our results can also be compared with \cite{SS}.  Finally,  several Hardy type inequalities and Heisenberg Pauli-Weyl uncertainty principles are derived as a consequence of the main result (see Subsection \ref{ssec32}).  However,  some basics and important results of conformable fractional calculus as required in this paper are presented in the next section (Section \ref{sec2}) as preliminaries.


 \section{Preliminaries}\label{sec2}
 This section gives basic definitions and properties which are fundamental to the concept of conformable fractional calculus as will be applied in the main results.  The notion of conformable divergence and Green's theorems which are vital to the proof of our results are also established.
 
\subsection{Basics concepts of conformable fractional calculus}
\begin{definition}\cite{Khal}
Given a function $u:[0,\infty)\to \mathbb{R}$. Then the conformable fractional derivative of order $\alpha$ is defined by
\begin{align*}
\left(T^\alpha_t u\right)(t) = \lim_{\delta\to 0} \frac{u(t+\delta t^{1-\alpha}) -u(t)}{\delta}
\end{align*}
for all $t>0$, $\alpha \in (0,1]$. If $u$ is $\alpha$-differentiable in some interval $(0,a)$, $a>0$, and $ \lim_{t\to 0}  \left(T^\alpha_t u\right)(t)$ exists, then it is defined as 
$$ (T^\alpha_t u)(0) =\lim_{t\to 0}  \left(T^\alpha_t u\right)(t).$$
\end{definition}
\noindent It has been well established that a function $u:[0,\infty)\to \mathbb{R}$ is continuous at $t_0>0$ if $u$ is $\alpha$-differentiable for $\alpha \in (0,1]$.

The next theorem collects those basic properties which conformable fractional derivative inherits from the classical derivative.
\begin{theorem}\cite{Abd,Khal}\label{thm23}
Let $\alpha \in (0,1]$, and let $f, g$ be $\alpha$-differentiable at a point $t>0$. Then 
\begin{enumerate}
\item $T^\alpha_t (a f + bg)(t) = a (T^\alpha_t f)(t) + b (T^\alpha_t g)(t), \ \ a, b\in \mathbb{R}.$
\item $T^\alpha_t (fg)(t) = g (T^\alpha_t f)(t) + f (T^\alpha_t g)(t).$
\item $\disp T^\alpha_t \left(\frac{f}{g}\right)(t) = \frac{g (T^\alpha_t f)(t) - f (T^\alpha_t g)(t)}{g^2(t)}, \ g\neq 0.$
\item $T^\alpha_t(C) = 0$ for all constant function $f(t)=C$.
\item $T^\alpha_t(t^s) = st^{s-\alpha}$ for $s\in \mathbb{R}$.
\item If in addition $f$ is differentiable, then 
$\disp (T^\alpha_t f)(t) = t^{1-\alpha}\frac{d u(t)}{dt}.$
\end{enumerate}
\end{theorem}

\begin{lemma}{\bf (Chain rule)}\cite{Abd,Khal}\label{lem24}
Assume $f$ is $\alpha$-differentiable with respect to $v$, and $v$ is $\alpha$-differentiable with respect to $x$. For $\alpha\in (0,1]$, we have 
$$T^\alpha_x\Big(f(v)\Big)(x)= \Big(T^\alpha_v f\Big)(v)\cdot v^{\alpha-1}\Big(T^\alpha_x v\Big)(x).$$
\end{lemma}

\begin{definition}{\bf ($\alpha$-fractional integral)} \cite{Abd,Khal}
Let $f$ be a continuous function on $[0,\infty)$, $t>a\ge 0$. Then for $\alpha \in (0,1]$, 
$$I^\alpha_a\Big(f(t)\Big) = I^1_a\Big(t^{\alpha -1}f(t)\Big)  =\int_a^t x^{\alpha-1}f(x)d x = \int_a^t f(x)d_\alpha x  .$$
\end{definition}
\noindent Here the integral is the usual Riemann improper integral. It easy to show that $T^\alpha_t  (I^\alpha_a f)(t) = f(t)$ whenever $f$ is continuous in the domain of $I^\alpha$:\\
$$T^\alpha_t  (I^\alpha_a f)(t) = t^{1-\alpha}\frac{d}{dt}I^\alpha_a(t) = t^{1-\alpha}\frac{d}{dt}\int_a^t x^{\alpha-1}f(x)d x = f(t).$$
Likewise,  $I^\alpha_a(T^\alpha_t  f)(t) = f(t)-f(a)$.

\begin{lemma}{\bf (Integration by Parts formula)}\cite{Abd,Khal}\label{lem25}
Suppose $f,g: [0,\infty)\to \mathbb{R}$ are $\alpha$-differentiable at a point $t>0$ for $\alpha\in (0,1]$. Then
$$\int_0^\infty \Big(T^\alpha_t f(t)\Big) g(t) d_\alpha t= f(t)g(t)\Big|_0^\infty  - \int_0^\infty f(t) \Big(T^\alpha_t g(t)\Big) d_\alpha t.$$
\end{lemma}

\begin{lemma}\cite{Ata}
Let $f$ be a non-constant differentiable function on an open interval. Then conformable derivative satisfies the following criteria
\begin{align*}
(a) \ \ \ \ \ T^{\alpha +\beta}_x\Big(f(x)\Big) & \neq T^\alpha_x \Big(T^\beta_x\Big(f(x)\Big)\Big) \ \ \ \text{for} \ \alpha, \beta \in (0,1). \\
(b) \ \ \ \ \ T^{\alpha +\beta}_x\Big(f(x)\Big) & \neq T^\beta_x \Big(T^\alpha_x\Big(f(x)\Big)\Big) \ \ \ \text{for} \ \alpha, \beta \in (0,1).\\
(c) \ \ \ \ \ T^{\alpha +\beta}_x\Big(f(x)\Big) & = T^\alpha_x \Big(T^\beta_x\Big(f(x)\Big)\Big) \ \ \ \text{for} \ \alpha, \in (0,1), \beta =1.
\end{align*}
\end{lemma}

Since many physical processes are modelled based on equations involving partial derivatives, it is inevitable to extend the above definition and properties to the case of differential of a function of several variables.

\begin{definition}\cite{Ata}
Let $f$ be a function of $n$-variables $x_1, x_2,\cdots ,x_n$. Then the conformable partial derivative of $f$ of order $\alpha\in (0,1]$ with respect to variable $x_k$, denoted by $D^\alpha_{x_k} := \frac{\partial^\alpha}{\partial x^\alpha_k}$, is defined as
\begin{align*}
D^\alpha_{x_k} f(\bar{x}) & =  \frac{\partial^\alpha f}{\partial  x^\alpha_k } (x_1, x_2,\cdots,x_n)\\
& =\lim_{\delta\to 0} \frac{f(x_1,\cdots,x_{k-1}, x_k+\delta x_k^{1-\alpha}, \cdots,x_n)- f(x_1,x_2, \cdots, x_n)}{\delta},
\end{align*}
$\bar{x}=(x_1,x_2,\cdots,x_n)$, $k=1,2,\cdots, n$.
\end{definition}
\noindent It has been noted that conformable partial derivative verifies the Clairaut criterion of mixed derivatives.

\begin{theorem}\cite{Ata}
Given a function $f(x,y)$ that is defined in the region $\Omega \subset \mathbb{R}^2$. Suppose $f$ has continuous conformable partial derivative of orders $\alpha$ and $\beta$, then
\begin{align*}
D^\alpha_{x}\Big(D^\beta_y\Big(f(x,y)\Big)\Big) = D^\beta_{y}\Big(D^\alpha_x\Big(f(x,y)\Big)\Big).
\end{align*}
\end{theorem}

 In the case that $f$ has conformable partial derivative of order $\alpha$ with respect to each variable $x_k$, $k=1,\cdots, n$. Then conformable vector can be defined at a point $q$ by 
$$D^\alpha_x f(q) = \Big(D^\alpha_{x_1}(f(q)),  D^\alpha_{x_2}(f(q)), \cdots, D^\alpha_{x_n}(f(q))\Big).$$

\noindent Consider the scalar field $f(\bar{x})$ and the vector field $\overrightarrow{F}(\bar{x})$ that are assumed to posses conformable partial derivative of order $\alpha$ with respect to all components $x_k$, $k=1,2,\cdots, n$.

\begin{definition}{\bf (Conformable gradient)}
The conformable gradient of order $\alpha$ as a vector field is given by
$$D^\alpha_x f(x) = \sum_{k=1}^n\Big(D^\alpha_{x_k} f) e_k,$$
where $e_k$ is the unit vector in the direction of $k$.
The conformable gradient of order $\alpha$ as a scalar field is given by
$$D^\alpha_x f(x) = \sum_{k=1}^n\Big(D^\alpha_{x_k} F_k).$$
\end{definition}

\begin{remark}
We note that conformable partial derivative (also conformable gradient) satisfies partial derivative versions of Theorem \ref{thm23}, Lemma \ref{lem24} and Lemma \ref{lem25}.
\end{remark}

\begin{definition}
By the above discussion,  anisotropic conformable fractional differential operator is therefore defined for continuous $\alpha$-differentiable function $f$ as
\begin{align*}
\sum^n_{k=1} D^\alpha_{x_k}\Big(| D^\alpha_{x_k} f(x)|^{p_k-2}D^\alpha_{x_k} f(x)\Big)  \ \ \text{for}\ \alpha \in (0,\infty],  \ \ p_k>1.
\end{align*}
\end{definition}

\noindent This can be written in the form of operator \eqref{e16}. We can now study fractional elliptic partial differential equations of the form

\begin{align*}
\sum^n_{k=1} D^\alpha_{x_k}\Big(| D^\alpha_{x_k} f(x)|^{p_k-2}D^\alpha_{x_k} f(x)\Big)  = g(x, f), \ \ \ x \in \Omega\subseteq\mathbb{R}^n
\end{align*}
in the appropriate fractional function spaces.

\subsection{Conformable Green's theorem}
Integration by parts formula, divergence theorem and Green's theorem within the framework of conformable fractional derivative will be applied severally. Then, there is a need to develop compatible divergence and Green's theorems for anisotropic conformable partial derivatives of order $\alpha$.

\begin{definition}\cite{Ata}
Let the vector field $F$ has the conformable partial derivatives of order $\beta$ on $\Omega \subseteq \mathbb{R}^n$. Then we denote by $P_x^\beta$ the vector
\begin{align*}
P_x^\beta F = \sum^n_{i=1} \Big\{e_{x_i}^T\Big(D^\beta_x(F)^T\Big)e_{x_i}  \Big\}e_{x_i}  = \sum^n_{i=1}\frac{\partial^\beta F_{x_i}}{\partial x^\beta_{x_i}} e_{x_i} .
\end{align*}
\end{definition}

\begin{definition}\cite{Ata}
Let the vector field $F$ has the conformable partial derivatives of order $\beta$ on an open region $\Omega$, $V\subseteq\Omega$ be simply connected and $S$ is the boundary surface of $V$ which is positively outward oriented. Then 
\begin{align*}
\iiint_{V} D^\alpha_x F \ d_\alpha V = \iint_S P^{\alpha-1}_x F \cdot n \ d_\alpha S.
\end{align*}
\end{definition}

\begin{remark}
This supports the fact that the conformable integral is anti-derivative of conformable derivative.
\end{remark}

\begin{lemma}{\bf (Conformable Green's Theorem)}\cite{Ata}
Let $C \subset \mathbb{R}^2$ be a simple positively oriented, piecewise smooth and closed region. Let $\Omega$ be the interior of $C$. If $f=f(x,y)$ and $g=g(x,y)$ have continuous conformable partial derivatives on $\Omega$. Then
\begin{align}\label{e21}
\iint_\Omega \Big( D^\alpha_x g - D^\alpha_y f\Big)d_\alpha S = \int_C D^{\alpha -1}_y fd_\alpha x + D^{\alpha-1}_x g d_\alpha y.
\end{align}
\end{lemma}

In what follows we consider a bounded open region $\Omega \subset \mathbb{R}^n$ with piecewise smooth and simple boundary. Note that the condition for the boundary to be simple amounts to $\partial\Omega$ being orientable. We say $\Omega \subset \mathbb{R}^n$ with this condition is said to be compatible.

\subsection{Green's identities}
Suppose $\Omega \subset \mathbb{R}^n$ is compatible, and $\alpha$-partial conformable fractional derivatives  $D^\alpha_{x_k}$ satisfy
\begin{align}\label{e22}
\sum^n_{k=1}\int_\Omega D^\alpha_{x_k} g_k(x) d_\alpha x  = \sum^n_{k=1}\int_{\partial\Omega} D^{\alpha-1}_{x_k}\Big(D^\alpha_{x_k} g_k(x)\Big)\cdot \nu  d_\alpha S 
\end{align}
for all $g_k \in \mathscr{D}^\alpha(\bar{\Omega})$, $k=1,2,\cdots,n$. Here $\bar{\Omega} =\Omega \cup \partial\Omega$,  $\mathscr{D}^\alpha(\bar{\Omega})$ denotes the space of all functions with continuous $\alpha$-partial conformable fractional derivative on $\Omega$ upto the boundary $\partial\Omega$, and $\nu$ is the outward pointing unit normal on $\partial\Omega$.
Next we prove Green's first and second identities for $\alpha$-partial conformable fractional derivative.

\begin{theorem}{\bf (Green's identities)}
Let $\Omega \subset \mathbb{R}^n$ be compatible,  we have
\begin{enumerate}
\item Green's first identity: Let $u, v \in \mathscr{D}^\alpha(\bar{\Omega})$, then 
\begin{align}
\int_\Omega\Big(D^\alpha_x u D^\alpha_x v + vD^\alpha_xD^\alpha_xu\Big) d_\alpha x  = \int_{\partial\Omega} v D^{\alpha-1}_x D^\alpha_x u \cdot \nu d_\alpha S.
\end{align}

\item Green's second identity: Let $u, v \in \mathscr{D}^\alpha(\bar{\Omega})$, then 
\begin{align}
\int_\Omega\Big(uD^\alpha_x u D^\alpha_x v - vD^\alpha_xD^\alpha_xu\Big) d_\alpha x  = \int_{\partial\Omega} \Big(u  D^{\alpha-1}_x D^\alpha_x v \cdot \nu - v D^{\alpha-1}_x D^\alpha_x u \cdot \nu \Big)d_\alpha S.
\end{align}
\end{enumerate}
\end{theorem}

\proof
Let $g_k = v D^\alpha_{x_k}u$, we have
\begin{align*}
\sum_{k=1}^n D^\alpha_{x_k} g_k & = \sum_{k=1}^n\Big(D^\alpha_{x_k}v D^\alpha_{x_k}u +v D^\alpha_{x_k}D^\alpha_{x_k}u \Big)\\
&= D^\alpha_{x}vD^\alpha_{x}u + vD^\alpha_{x}D^\alpha_{x}u.
\end{align*}
Applying the divergence formula \eqref{e22} leads to
\begin{align*}
\int_\Omega\Big(D^\alpha_{x}vD^\alpha_{x}u + vD^\alpha_{x}D^\alpha_{x}u \Big)d_\alpha x  & =  \sum_{k=1}^n D^\alpha_{x_k} g_k\\
& = \sum^n_{k=1}\int_{\partial\Omega} D^{\alpha-1}_{x_k}\Big(D^\alpha_{x_k} g_k(x)\Big)\cdot \nu  d_\alpha S\\
&= \int_{\partial\Omega} v D^{\alpha-1}_{x} D^{\alpha}_{x} u\cdot \nu d_\alpha S.
\end{align*}
Applying the Green's first identity by swapping the positions of $u$ and $v$ we have

\begin{align}
\int_\Omega\Big(D^\alpha_x u D^\alpha_x v + vD^\alpha_xD^\alpha_xu\Big) d_\alpha x  & = \int_{\partial\Omega} v D^{\alpha-1}_x D^\alpha_x u \cdot \nu d_\alpha S\\
\int_\Omega\Big(D^\alpha_x v D^\alpha_x u + uD^\alpha_xD^\alpha_xv\Big) d_\alpha x  & = \int_{\partial\Omega} v D^{\alpha-1}_x D^\alpha_x u \cdot \nu d_\alpha S.
\end{align}
Subtracting one of the equations from the other yields the desire result.
\qed

\begin{remark}
If $v=1$ in these Green's identity we obtain the following analogue of Gauss mean value formula for $\alpha$-conformable harmonic function satisfying $D^\alpha_x D^\alpha_xu =0$ in a compatible domain:
$$\int_{\partial\Omega}  D^{\alpha-1}_x \Big(D^\alpha_x u\Big) \cdot \nu d_\alpha S = 0.$$
\end{remark}

\section{Anisotropic conformable Hardy type inequalities with weights}
First, we prove some Picone type identities for anisotropic fractional gradient operator on a compatible domain.  The generalized conformable fractional Hardy inequalities are derived as a consequence of the anisotropic Picone identities.  Comparing our results with those recently obtained for single valued functions in \cite{GSAZ1, GSAZ2}, it can be said that the approach adopted in this paper appears more general as  multidimensional functions are considered and it is able to take care of the boundary term. 

\subsection{Anisotropic Picone type identity  and Hardy type inequalities}
\begin{proposition}{\bf (Anisotropic conformable Picone type identity)}\label{prop31}
Let $u\ge 0$ and $v>0$ be $\alpha$-order conformable differentiable functions a.e. in an open domain $\Omega\subset \mathbb{R}^n$. Define
\begin{align*}
\mathcal{R}(u,v) = \sum_{k=1}^n |D^\alpha_{x_k} u|^{p_k} -  \sum_{k=1}^n D^\alpha_{x_k}\left( \frac{u^{p_k}}{v^{p_k-1}}\right)|D^\alpha_{x_k} v|^{p_k-2}D^\alpha_{x_k} v,
\end{align*} 
\begin{align*}
\mathcal{L}(u,v) = \sum_{k=1}^n |D^\alpha_{x_k} u|^{p_k} 
+ \sum_{k=1}^n (p_k-1)\frac{u^{p_k}}{v^{p_k}} |D^\alpha_{x_k} v|^{p_k}
-  \sum_{k=1}^n  p_k \frac{u^{p_k-1}}{v^{p_k-1}} |D^\alpha_{x_k} v|^{p_k-2} D^\alpha_{x_k} u D^\alpha_{x_k} v,
\end{align*} 
where $0<\alpha\le 1$ and $p_k>1$, $k=1,2,\cdots, n$. Then
$\mathcal{R}(u,v) =\mathcal{L}(u,v)\ge 0$. Moreover, $\mathcal{L}(u,v) =0$ a.e. in $\Omega$ if and only if $u/v=c$ for a positive constant $c$.
\end{proposition}

\proof
By the chain rule for conformable partial derivative we compute
\begin{align*}
D^\alpha_{x_k} (u^{p_k}) &= D^\alpha_u(u^{p_k})u^{\alpha-1}D^\alpha_{x_k} (u) = p_ku^{p_k-1} D^\alpha_{x_k} u.\\
D^\alpha_{x_k} (v^{p_k-1}) & = D^\alpha_v(v^{p_k-1})v^{\alpha-1}D^\alpha_{x_k} (v) = (p_k-1)v^{p_k-2} D^\alpha_{x_k} v.
\end{align*}
Now evaluating by quotient rule for conformable partial derivative gives
\begin{align*}
D^\alpha_{x_k} \left(\frac{u^{p_k}}{v^{p_k-1}}\right) = \frac{p_ku^{p_k-1} D^\alpha_{x_k} u}{v^{p_k-1}} - \frac{(p_k-1)u^{p_k} D^\alpha_{x_k} v}{v^{p_k}}.
\end{align*}
Substituting this into the expression for $\mathcal{R}(u,v)$ we arrive at $\mathcal{R}(u,v) =\mathcal{L}(u,v)$.
To prove that $\mathcal{L}(u,v)\ge 0$, the expression for $\mathcal{L}(u,v)$ can be broken down as follows
\begin{align*}
\mathcal{L}(u,v) = \mathcal{A}_1(u,v) + \mathcal{A}_2(u,v),
\end{align*}
where
\begin{align*}
\mathcal{A}_1(u,v) & = \sum_{k=1}^n |D^\alpha_{x_k} u|^{p_k} - \sum_{k=1}^n  p_k \frac{u^{p_k-1}}{v^{p_k-1}} |D^\alpha_{x_k} v|^{p_k-1} |D^\alpha_{x_k} u|  + \sum_{k=1}^n (p_k-1)\frac{u^{p_k}}{v^{p_k}} |D^\alpha_{x_k} v|^{p_k},  \\ \ \\
\mathcal{A}_2(u,v) & =  \sum_{k=1}^n  p_k \frac{u^{p_k-1}}{v^{p_k-1}} |D^\alpha_{x_k} v|^{p_k-2} \Big\{ |D^\alpha_{x_k} v| |D^\alpha_{x_k} u| -D^\alpha_{x_k} u D^\alpha_{x_k} v\Big\}.
\end{align*}

Recall the Young's inequality for real numbers $a,b\ge0$ and exponents $p_k>1$, $q_k>1$ satisfying $1/q_k = (p_k-1)/p_k$:
\begin{align}\label{e31}
ab \le (1/p_k)a^{p_k} +((p_k-1)/p_k)b^{p_k/(p_k-1)}
\end{align}
with equality if and only if $a=b^{1/(p_k-1)}$ for all $k=1,2,\cdots, n$.

 We claim that $\mathcal{A}_1(u,v)\ge 0$ and verify it by the application of the Young's inequality above. Set $a = |D^\alpha_{x_k} u|$ and $b = \left(\frac{u}{v}|D^\alpha_{x_k} v|\right)^{p_k-1}$ and then by \eqref{e31} we have 
 \begin{align*}
  \frac{u^{p_k-1}}{v^{p_k-1}} |D^\alpha_{x_k} v|^{p_k-1} |D^\alpha_{x_k} u| \le \frac{1}{p_k} |D^\alpha_{x_k} u|^{p_k} + \frac{p_k-1}{p_k}\left(\frac{u}{v} |D^\alpha_{x_k} v|\right)^{p_k}.
 \end{align*}
 Therefore 
 \begin{align*}
 \mathcal{A}_1(u,v) & = \sum_{k=1}^n p_k\left[ \frac{1}{p_k} |D^\alpha_{x_k} u|^{p_k} + \frac{p_k-1}{p_k}\left(\frac{u}{v} |D^\alpha_{x_k} v|\right)^{p_k}\right] - \sum_{k=1}^n  p_k \frac{u^{p_k-1}}{v^{p_k-1}} |D^\alpha_{x_k} v|^{p_k-1} |D^\alpha_{x_k} u|\\
 & \ge 0.
\end{align*} 
 We also see that $ \mathcal{A}_2(u,v)\ge 0$ by referring to Cauchy-Schwarz  inequality in the form 
 $$D^\alpha_{x_k} u D^\alpha_{x_k} v \le |D^\alpha_{x_k} u| |D^\alpha_{x_k} v| .$$
 It can therefore be concluded that 
 $\mathcal{L}(u,v) = \mathcal{A}_1(u,v) + \mathcal{A}_2(u,v) \ge 0$.

It is straightforward to see that $u=cv$ yields $\mathcal{R}(u,v)=0$. Now, we must prove that $\mathcal{L}(u,v)=0$ if and only if $u=cv$. From the above analysis we observe that $\mathcal{L}(u,v)=0$ if and only if 
\begin{align}\label{e32}
|D^\alpha_{x_k} u| =\frac{u}{v} |D^\alpha_{x_k} v|, \ \ \ \ \ k=1,2,\cdots, n
\end{align}
\begin{align}\label{e33}
|D^\alpha_{x_k} u|  |D^\alpha_{x_k} v| = D^\alpha_{x_k} u   D^\alpha_{x_k} v, \ \ \ \ \ k=1,2,\cdots, n.
\end{align}
Combining the last two identities yields $D^\alpha_{x_k} u/ D^\alpha_{x_k} v = u/v  = c$ for $c>0$, $k=1,2,\cdots, n$, which implies $D^\alpha_{x_k} (u/v)=0$ .

Since $u(x)\ge 0$ and $\mathcal{L}(u,v)(x_0)=0$, $x_0\in \Omega$ we need to verify two cases, namely, $u(x_0)>0$ and $u(x_0)=0$:

\noindent (a).  If $u(x_0)>0$, then $\mathcal{L}(u,v)=0$ for all $x_0\in \Omega$ and we conclude that \eqref{e32} and \eqref{e33} hold which when combined give $u=cv$ a.e. for constant $c>0$ for all $x_0\in\Omega$.

\noindent (b). If $u(x_0)=0$, define the set $\Omega^*=\{x\in\Omega: u(x_0)=0\}$, and suppose $\Omega^*\ne \Omega$. Then there exists a sequence $x_m \notin \Omega^*$ such that $x_m\to x_0$. In particular, $u(x_m)\ne 0$, and hence by the case (a) $u(x_m) =cv(x_m)$. Passing to the limit we get $u(x_0)=cv(x_0)$ which implies $c=0$ since $u(x_0)=0$ and $v(x_0)>0$.  By the case (a) again we know that $u=cv$ and $u=0$ for all $x\in \Omega\setminus\Omega^*$, then it is impossible to have $c=0$. This contradiction implies that $\Omega^*=\Omega$.
This completes the proof.

\qed

The stage is now set to present and prove the weighted fractional Hardy type inequality which is our main result.

\begin{theorem}{\bf (Anisotropic conformable Hardy type inequality)}\label{thm32}
Let $W_k(x)\ge 0$ and $H_k(x)\ge 0$ be two weight functions, where $k=1,2,\cdots n$, such that a continuous $\alpha$-differentiable function $v>0$ a.e. in a compatible domain $\Omega \subset \mathbb{R}^n$ (that is $0<v\in \mathscr{D}^\alpha(\Omega)\cap C(\bar{\Omega})$) satisfies conformable partial differential inequality
\begin{align}\label{e34}
- D^\alpha_{x_k} \Big(W_k(x)|D^\alpha_{x_k}v|^{p_k-2}D^\alpha_{x_k} v\Big) \ge L_k H_k(x) v^{p_k-1}
\end{align}
for $k=1,2,\cdots n$. Then for all functions $u\ge 0$ which are continuous $\alpha$-differentiable in $\Omega$, we have 
\begin{align}\label{e35}
 \sum_{k=1}^n \int_\Omega W_k(x) & |D^\alpha_{x_k} u|^{p_k} d_\alpha x \ge   \sum_{k=1}^n L_k \int_\Omega H_k(x) |u|^{p_k}d_\alpha x  \nonumber \\
 & +  \sum_{k=1}^n  \int_{\partial\Omega} \frac{u^{p_k}}{v^{p_k-1}}  D_{x_k}^{\alpha-1} D_{x_k}^{\alpha} \Big(W_k(x)|D^\alpha_{x_k}v|^{p_k-2}D^\alpha_{x_k} v\Big) \cdot \nu d_\alpha A.
\end{align}
\end{theorem}

If $u$ vanishes on the boundary $\partial\Omega$, then the last term on the right hand side of  \eqref{e35} will vanish and we then have
\begin{align}\label{e36}
 \sum_{k=1}^n \int_\Omega W_k(x) & |D^\alpha_{x_k} u|^{p_k} d_\alpha x \ge   \sum_{k=1}^n L_k \int_\Omega H_k(x) |u|^{p_k}d_\alpha x.
\end{align}

\proof
In this proof we apply the anisotropic Picone identities in Proposition \ref{prop31}, Green's identity and conformable partial differential inequality \eqref{e35} as follows:
\begin{align*}
0 &\le \sum_{k=1}^n \int_\Omega W_k(x)\mathcal{L}(u,v) d_\alpha x \\
 & =  \sum_{k=1}^n \int_\Omega W_k(x) |D^\alpha_{x_k} u|^{p_k}  d_\alpha x -  \sum_{k=1}^n \int_\Omega W_k(x) D^\alpha_{x_k}\left( \frac{u^{p_k}}{v^{p_k-1}}\right)|D^\alpha_{x_k} v|^{p_k-2}D^\alpha_{x_k} v d_\alpha x\\
 & = \sum_{k=1}^n \int_\Omega W_k(x) |D^\alpha_{x_k} u|^{p_k}  d_\alpha x +  \sum_{k=1}^n \int_\Omega  \frac{u^{p_k}}{v^{p_k-1}} D^\alpha_{x_k}\Big(W_k(x) |D^\alpha_{x_k} v|^{p_k-2}D^\alpha_{x_k} v\Big) d_\alpha x\\
 & \hspace{2cm} +  \sum_{k=1}^n \int_\Omega  \frac{u^{p_k}}{v^{p_k-1}} D^{\alpha-1}_{x_k}\ D^\alpha_{x_k}\Big(W_k(x) |D^\alpha_{x_k} v|^{p_k-2}D^\alpha_{x_k} v\Big) \cdot \nu d_\alpha A\\
 & \le  \sum_{k=1}^n \int_\Omega W_k(x) |D^\alpha_{x_k} u|^{p_k}  d_\alpha x -   \sum_{k=1}^n L_k \int_\Omega   H_k(x) u^{p_k} d_\alpha x\\
 & \hspace{2cm} +  \sum_{k=1}^n \int_\Omega  \frac{u^{p_k}}{v^{p_k-1}} D^{\alpha-1}_{x_k}\ D^\alpha_{x_k}\Big(W_k(x) |D^\alpha_{x_k} v|^{p_k-2}D^\alpha_{x_k} v\Big) \cdot \nu d_\alpha A.
\end{align*}
This proves the required statement.

\qed

\subsection{Consequences of the conformable  weighted  Hardy type inequalities}\label{ssec32}

\begin{theorem}\label{thm33}
Let $\Omega\subset\mathbb{R}^n$ be a compatible domain (or open bounded domain). Let $m\in \mathbb{R}$, $a>\alpha\in (0,1]$, $1<p_k<a+m$ and $p_k(1-\alpha)\ge a-\alpha$ for $k=1,2,\cdots n$.  Then for all  $\alpha$-partial conformable differentiable function  $u\in \mathscr{D}^\alpha(\Omega\setminus \{x_k=0\}, k=1,2,\cdots n)$ we have
\begin{align}\label{e37}
 \sum_{k=1}^n \int_\Omega |x_k|^{m}  |D^\alpha_{x_k} u|^{p_k} d_\alpha x \ge   \sum_{k=1}^n \left(\frac{m+a-p_k}{p_k}\right)^{p_k} \int_\Omega \frac{|u|^{p_k}}{|x_k|^{\alpha p_k-m}}d_\alpha x.
\end{align}
\end{theorem}

\proof
Without loss of generality we consider $u\ge 0$, $\alpha$-partial conformable differentiable. Define an auxiliary function
$$v=\prod_{k=1}^n |x_k|^{\beta_k} = |x_j|^{\beta_j}V_k,$$
where $V_k=\prod_{k=1,k\neq j}^n |x_k|^{\beta_k}$ and $\beta_k = -(m+a-p_k)/p_k$. Also define $W_k(x) =|x_k|^m, m\in \mathbb{R}$.
Then by a straightforward computation we have 
\begin{align*}
D^\alpha_{x_k} v & = \beta_k V_k|x_k|^{\beta_k-\alpha}\\
|D^\alpha_{x_k} v|^{p_k-2} & =|\beta_k|^{p_k-2}V_k^{p_k-2}|x_k|^{\beta_kp_k-2\beta_k-\alpha p_k+2\alpha}\\
|D^\alpha_{x_k} v|^{p_k-2}D^\alpha_{x_k} v &= |\beta_k|^{p_k-2}\beta_k V_k^{p_k-1} |x_k|^{\beta_kp_k- \beta_k-\alpha p_k+\alpha}.
\end{align*}
We compute
\begin{align*}
D^\alpha_{x_k}& \Big( W_k(x) |D^\alpha_{x_k} v|^{p_k-2}D^\alpha_{x_k} v \Big)\\
 & = |\beta_k|^{p_k-2}\beta_k D^\alpha_{x_k}\Big(V_k^{p_k-1} |x_k|^{\beta_kp_k- \beta_k-\alpha p_k+\alpha+m}\Big)\\
&=  |\beta_k|^{p_k-2}\beta_k\Big(-\beta_k+p_k(\beta_k-\alpha)+\alpha+m\Big) V_k^{p_k-1}|x_k|^{\beta_k(p_k-1)-\alpha p_k+m}\\
& =-  |\beta_k|^{p_k}  (V_k|x_k|^{\beta_k})^{p_k-1}|x_k|^{-\alpha p_k+m}\\
& \hspace{1cm} +  |\beta_k|^{p_k-2}\beta_k\Big(p_k(\beta_k-\alpha)+\alpha+m\Big) (V_k|x_k|^{\beta_k})^{p_k-1}|x_k|^{-\alpha p_k+m}.
\end{align*}
Now using $\beta_k = -(m+a-p_k)/p_k$  and $1<p_k<a+m$ we get 
\begin{align*}
-D^\alpha_{x_k}& \Big( W_k(x) |D^\alpha_{x_k} v|^{p_k-2}D^\alpha_{x_k} v \Big) = \left|\frac{m+a-p_k}{p_k}\right|^{p_k} \frac{v^{p_k-1}}{|x_k|^{\alpha p_k-m}} \\
&+ \left|\frac{m+a-p_k}{p_k}\right|^{p_k-2} \left(\frac{m+a-p_k}{p_k}\right)\Big(p_k(\beta_k-\alpha)+\alpha+m\Big)\frac{v^{p_k-1}}{|x_k|^{\alpha p_k-m}} \\
& \ge  \left|\frac{m+a-p_k}{p_k}\right|^{p_k} \frac{v^{p_k-1}}{|x_k|^{\alpha p_k-m}},
\end{align*}
where the inequality is due to the fact that $p_k(\beta_k-\alpha)+\alpha+m>0$ which comes as a result of the condition $p_k(1-\alpha)\ge a-\alpha$.  This has clearly fulfilled the condition \eqref{e34}. Then plugging in the following data
$$W_k(x) =|x_k|^m, \ \ L_k =  \left|\frac{m+a-p_k}{p_k}\right|^{p_k}\ \ \text{and} \ \ H_k(x)=  \frac{1}{|x_k|^{\alpha p_k-m}}$$
into \eqref{e35} we arrived at the desired inequality.

\qed

The next two corollaries give some special cases of \eqref{e37}:

\begin{corollary}
With the conditions of Theorem \ref{thm33} the following inequalities hold: 
\begin{align}\label{e38}
 \sum_{k=1}^n \int_\Omega   |D^\alpha_{x_k} u|^{p_k} d_\alpha x \ge   \sum_{k=1}^n \left|\frac{a-p_k}{p_k}\right|^{p_k} \int_\Omega \frac{|u|^{p_k}}{|x_k|^{\alpha p_k}}d_\alpha x
\end{align}
and
\begin{align}\label{e39}
 \sum_{k=1}^n \int_\Omega   |D^\alpha_{x_k} u|^{p_k} d_\alpha x \ge   \sum_{k=1}^n \left(\frac{\alpha(p_k-1)}{p_k}\right)^{p_k} \int_\Omega \frac{|u|^{p_k}}{|x_k|^{\alpha p_k}}d_\alpha x.
\end{align}
\end{corollary}

\proof
Setting $m=0$ $\implies$  $W_k(x)=1$ in \eqref{e37} and gives rise to \eqref{e38}. Furthermore, using the condition $p_k(1-\alpha)\ge a-\alpha$, one sees that $|(a-p_k)/p_k| \ge [\alpha(p_k-1)/p_k]$ which proves \eqref{e39}.

\qed

\begin{corollary}
With the conditions of Theorem \ref{thm33} the following inequality holds: 
\begin{align}\label{e310}
 \sum_{k=1}^n \int_\Omega |x_k|^{\alpha}  |D^\alpha_{x_k} u|^{p_k} d_\alpha x \ge   \sum_{k=1}^n \left(\frac{\alpha+a-p_k}{p_k}\right)^{p_k} \int_\Omega \frac{|u|^{p_k}}{|x_k|^{\alpha( p_k-1)}}d_\alpha x.
\end{align}
\end{corollary}

\proof
Set $m=\alpha$ $\implies$ $W_k(x)=|x|^\alpha$, $\alpha\in (0,1]$ in  \eqref{e37}. 

\qed

\begin{theorem}\label{thm36}
Let $\Omega\subset\mathbb{R}^n$ be a compatible domain. Then the following inequality holds
\begin{align}\label{e311}
 \sum_{k=1}^n \int_\Omega  |D^\alpha_{x_k} u|^{p_k} d_\alpha x \ge   \sum_{k=1}^n (p_k-1)|m|^{p_k} \int_\Omega \frac{|u|^{p_k}}{|x_k|^{(\alpha-1) p_k}}d_\alpha x
\end{align}
for all $\alpha$-partial conformable differentiable functions $u\in \mathscr{D}_0^\alpha(\Omega\setminus \{x_k=0\}, k=1,2,\cdots n)$,  where $m<0$.
\end{theorem}

\noindent This theorem can be proved using similar approach as used in the proof of Theorem \ref{thm33}.
\proof
Here we consider the weight function $W_k(x)=1$ and the auxiliary function $v=e^{m|x_k|}$,  $m<0$, and compute as follows:
\begin{align*}
 D^\alpha_{x_k} v &= m|x_k|^{1-\alpha}e^{m|x_k|} \\
 |D^\alpha_{x_k} v|^{p_k-2}  &=|m|^{p_k-2}|x_k|^{(1-\alpha)(p_k-2)} e^{m(p_k-2)|x_k|}\\
 |D^\alpha_{x_k} v|^{p_k-2}D^\alpha_{x_k} v &=|m|^{p_k-2}m |x_k|^{(1-\alpha)(p_k-1)} e^{m(p_k-1)|x_k|}
\end{align*}
and then 
\begin{align*}
 D^\alpha_{x_k}(W_k |D^\alpha_{x_k} v|^{p_k-2}D^\alpha_{x_k} v) &=   |m|^{p_k-2}m (1-\alpha)(p_k-1) |x_k|^{(1-\alpha)(p_k-1)-\alpha} e^{m(p_k-1)|x_k|}\\
 & -|m|^{p_k}(p_k-1)|x_k|^{(1-\alpha)p_k} e^{m(p_k-1)|x_k|}.
\end{align*} 
Therefore 
\begin{align*}
- D^\alpha_{x_k}(W_k |D^\alpha_{x_k} v|^{p_k-2}D^\alpha_{x_k} v) \ge (p_k-1)|m|^{p_k} \frac{v^{p_k-1}}{|x_k|^{p_k(\alpha-1)}}
\end{align*}
which has satisfied the hypothesis of Theorem \ref{thm32}. So we choose
$$W_k =1, \ \ L_k=(p_k-1)|m|^{p_k} \ \ \text{and}\ \ H_k(x)= \frac{1}{|x_k|^{p_k(\alpha-1)}}$$
and then by \eqref{e36} the desired inequality is obtained. 

\qed

\begin{remark}\label{rem37}
Choosing $m =-\frac{n-p_k}{p_k}$ with $1<p_k<n$ gives an interesting anisotropic conformable Hardy inequality with weight. If in addition one sets $p_k=2$, $n\ge 3$, \eqref{e310} results to an isotropic version:
\begin{align}\label{e312}
  \int_\Omega  |D^\alpha_{x} u|^2 d_\alpha x \ge   \left(\frac{n-2}{2}\right)^2 \int_\Omega \frac{u^2}{|x|^{2(\alpha-1)}}d_\alpha x
\end{align}
which is well known in the classical setting with $\alpha=1$.
\end{remark}

\subsection{Heisenberg-Pauli-Weyl (HPW) uncertainty principles}
Lastly, it is demonstrated here that Heisenberg-Pauli-Weyl (HPW) uncertainty principles can be derived as an immediate consequence of the above conformable fractional Hardy inequality. The classical uncertainty principle of quantum mechanics says that some pairs of physical quantities, such as position and momentum of a particle,  cannot be determined exactly at the same time but only with an 'uncertainty'.  As a theorem in Euclidean harmonic analysis,  uncertainty principle expresses impossibility of simultaneously smallness of a nonzero function $f$ and its Fourier transform $\hat{f}$ (where $\hat{f}(y)= (2\pi)^{-n/2} \int_{\mathbb{R}^n} f(x) e^{i\langle x, y\rangle} dx$ and $\|f\|_2=1 = \|\hat{f}\|_2$).

\begin{corollary}{\bf (Fractional HPW uncertainty inequality)}\label{cor38}
Let $\Omega\subset\mathbb{R}^n$ , $n>2$, be a compatible domain. Then for all $\alpha$-partial conformable differentiable functions $u\in \mathscr{D}_0^\alpha(\Omega\setminus \{x=0\})$,  $\alpha \in (0,1]$, the following inequality  
\begin{align}\label{e313}
 \left( \int_\Omega|x|^{2(\alpha -1)} |u|^2  d_\alpha x \right)\left( \int_\Omega  |D^\alpha_{x} u|^2 d_\alpha x\right)  \ge   \left(\frac{n-2}{2}\right)^2  \left(\int_\Omega  |u|^2d_\alpha x\right)^2
\end{align}
holds.
\end{corollary}

\proof
Starting with Cauchy-Schwarz inequality we have
\begin{align*}
\int_\Omega |u|^2 d_\alpha x \le \left(\int_\Omega \frac{|u|^2}{|x|^{2(\alpha-1)}}d_\alpha x\right)^{1/2} \left(\int_\Omega  |x|^{2(\alpha-1)} |u|^2d_\alpha x\right)^{1/2}.
\end{align*}
Applying the conformable Hardy inequality \eqref{e311} we obtain
\begin{align*}
\int_\Omega |u|^2 d_\alpha x \le \frac{2}{n-2} \left(\int_\Omega  |D^\alpha_x u|^2d_\alpha x\right)^{1/2} \left(\int_\Omega  |x|^{2(\alpha-1)} |u|^2d_\alpha x\right)^{1/2}
\end{align*}
which is the desired inequality.

\qed

\begin{theorem}{\bf (Anisotropic fractional HPW uncertainty inequality)}\label{thm39}
Let $\Omega\subset\mathbb{R}^n$  be a compatible domain, $1<p_k,q_k<\infty$ and $1/p_k+1/q_k=1$ for all $k=1,2,\cdots, n$. Then for all $\alpha$-partial conformable differentiable functions $u\in \mathscr{D}_0^\alpha(\Omega\setminus \{x_k=0\}, k=1,2,\cdots n)$,  $\alpha \in (0,1]$, the following inequality  
\begin{align}\label{e314}
\sum_{k=1}^n  \left(\int_\Omega|x|^{q_k(\alpha -1)} |u|^{p_k}  d_\alpha x \right)^{\frac{1}{q_k}}\left(\int_\Omega  |D^\alpha_{x} u|^{p_k} d_\alpha x\right)^{\frac{1}{p_k}}  \ge    \sum_{k=1}^n (p_k-1)|m| \int_\Omega |u|^{p_k} d_\alpha x
\end{align}
holds.
\end{theorem}

\proof The proof is similar to that of Corollary \ref{cor38} by applying H\"older inequality with $1/p_k+1/q_k=1$, and then Theorem \ref{thm36}.

\qed

If $m$ is chosen as in Remark \ref{rem37} and $p_k=2$, then  Corollary \ref{cor38} becomes a special case of Theorem \ref{thm39}.

 \section*{Declarations}

 \subsection*{Ethics approval and consent to participate}
Not applicable.

 \subsection*{Competing  interests}
 The author declares that there is no competing interests.
 
\subsection*{Author's contributions}
AA wrote the manuscript and YOA edited the manuscript. 
Both authors approved the manuscript.

\subsection*{Funding}
This work  does not receive any funding.

\subsection*{Availability of data and material}
This manuscript has no associated data.

\subsection*{Acknowledgement} 
The authors thank their respective institutions. 


\end{document}